\documentclass[12pt]{article}
\usepackage{amsfonts,epsfig,subfigure,fancybox}
\def\qed{$\qquad \Box$}
\def\hh{{\hbox{I{\kern -0.22em}H}}}

\def\op{{\cal L}}
\def\wdh{\widehat}

\def\tr{\hbox{tr}}
\newcommand{\diag}{{\rm diag}}

\newcommand{\nd}{\noindent}
\newcommand{\wdt}{\widetilde}
\def\para#1{\vskip 0.4\baselineskip\noindent{\bf #1}}

\newcommand{\disp}{\displaystyle}
\newcommand{\ind}{\mbox{1}\kern-.25em \mbox{I}}

\makeatletter
\@addtoreset{equation}{section}
\renewcommand{\thesection}{\arabic{section}}
\renewcommand{\theequation}{\thesection.\arabic{equation}}
\newcommand{\beq}[1]{\begin{equation} \label{#1}}
\newcommand{\eeq}{\end{equation}}
\newcommand{\bed}{\begin{displaymath}}
\newcommand{\eed}{\end{displaymath}}
\newcommand{\bea}{\bed\begin{array}{rl}}
\newcommand{\eea}{\end{array}\eed}
\newcommand{\ad}{&\!\!\!\disp}
\newcommand{\aad}{&\disp}
\newcommand{\barray}{\begin{array}{ll}}
\newcommand{\earray}{\end{array}}

\newcommand{\ph}{\varphi}

\newcommand{\al}{\alpha}

\newcommand{\ga}{\gamma}

\def\k0{\kappa_0}

\newcommand{\rr}{{\Bbb R}}
\newcommand{\M}{{\mathcal M}}
\newcommand{\A}{{\mathcal A}}
\newcommand{\dE}{{\mathbb E}}
\newcommand{\dP}{{\mathbb P}}
\newcommand{\F}{{\mathcal F}}

\newcommand{\lbar}{\overline}

\def\one{{\hbox{1{\kern -0.35em}1}}}

\def\lan{\big\langle}
\def\ran{\big\rangle}

\newcommand{\bdd}{\hspace*{-0.08in}{\bf.}\hspace*{0.05in}}
\newtheorem{thm}{Theorem}[section]
\newtheorem{prop}[thm]{Proposition}
\newtheorem{lem}[thm]{Lemma}
\newtheorem{rem}[thm]{Remark}
\newtheorem{cor}[thm]{Corollary}
\newtheorem{defn}[thm]{Definition}

\newcommand{\thmref}[1]{Theorem~{\rm \ref{#1}}}
\newcommand{\lemref}[1]{Lemma~{\rm \ref{#1}}}

\newcommand{\propref}[1]{Proposition~{\rm \ref{#1}}}

\def\({\left(}
\def\){\right)}

\def\l{\left|}
\def\r{\right|}

\begin{document}

\title{Almost Sure Stabilization for Adaptive
Controls of Regime-switching LQ
Systems with A Hidden Markov
Chain}
\date{}
\author{Bernard Bercu\thanks{Universit\'e Bordeaux 1,
Institut de Math\'ematiques de Bordeaux,
UMR 5251 and INRIA Bordeaux Sud-Ouest, Team CQFD,
351 cours de la Lib\'eration, 33405 Talence cedex, France.
Email: {\tt Bernard.Bercu@math.u-bordeaux1.fr}}
\and Fran\c{c}ois Dufour\thanks{Universit\'e Bordeaux 1,
Institut de Math\'ematiques de Bordeaux,
UMR 5251 and INRIA Bordeaux Sud-Ouest, Team CQFD,
351 cours de la Lib\'eration, 33405 Talence cedex, France.
Email:
{\tt Francois.Dufour@math.u-bordeaux1.fr}}
\and G. George Yin\thanks{Department of Mathematics, Wayne State
University, Detroit, Michigan 48202. Email: {\tt
gyin@math.wayne.edu}}}

\maketitle
\begin{abstract}
This work is devoted to the almost sure stabilization
of adaptive control systems that involve an unknown Markov chain.
The control system displays continuous dynamics
represented by differential equations and discrete events given by a hidden
Markov chain. Different from previous work on stabilization of
adaptive controlled systems with a hidden Markov chain, where
average criteria were considered,
this work focuses on the almost
sure stabilization or sample path stabilization of the underlying
processes. Under simple conditions, it is shown that as long as
the feedback controls have linear growth in the continuous
component, the resulting process is regular. Moreover, by appropriate
choice of the Lyapunov functions, it is shown that the adaptive
system is stabilizable almost surely. As a by-product, it is also
established that the controlled process is positive recurrent.

\vskip 0.2 true in
\nd{\bf Key words.}
Adaptive control, hidden Markov chain, almost sure stabilization.

\vskip 0.2 true in
\nd{\bf Abbreviated title.}  Almost Sure
Stabilization of Adaptive Controls with Switching

\end{abstract}

\newpage
\setlength{\baselineskip}{0.28in}
\section{Introduction}
This work deals with almost sure stabilization of
adaptive control systems in continuous-time
with an unknown parameter process that is a hidden Markov chain.
The systems belong to the class of partially observed control
systems. Naturally, one estimates the parameter process
by using nonlinear filtering techniques and then uses the estimator
in the systems in order to design adaptive control strategies.
The motivation of our study stems from consideration of the following
problem. Let us begin with a hybrid linear quadratic (LQ) problem
$$\dot X(t) = A_{\al(t)}X(t)+ B_{\al(t)} U(t)$$
where $\al(t)$ is a continuous-time Markov chain taking values in a
finite set $\M=\{1,\ldots,m\}$, $A_i $ and $B_i$ for $i\in \M$ are
matrices with compatible dimensions, and $U(t)$ is the control process.
One can observe that different from the traditional
setup of LQ problems, the system matrices $A_i$ and $B_i$ are both
subject to random switching influence. At any given instance,
these coefficient matrices are chosen from a set $\M$ with
a finite number of candidates. The selection rule is dictated by the
modulating switching process $\al(t)$ that jump changes from one
state to another at random times. Such systems have enjoyed
numerous applications in emerging application areas as financial
engineering, wireless communications, as well as in existing
applications. A particular important problem concerns the
asymptotic behavior of such systems when they are in operations
for a long time. Our interest lies in finding
admissible controls so that the resulting system will be almost surely
stabilized.
An added difficulty is that the process $X(t)$ can only be observed with an
additive noise
$$dX(t)=[A_{\al(t)}X(t)+ B_{\al(t)}U(t)]dt +dW(t).$$
For such partially observed systems, it is natural to use nonlinear filtering
techniques. The associated filter is known as the Wonham filter
\cite{Wonham65}, which is one of a
handful of finite dimensional filters in
existence.

Linear quadratic (LQ) regulators appear to present rather simple
structures. Meanwhile, there are so many applications that can be
described by such processes. We refer the reader to
\cite{BlairS,CainesC,JiC,Mariton} for some recent work on the
associated control, estimation, and optimization problems for
hybrid systems. Emerging applications have also been found in
manufacturing systems, in which a Markov chain is used to
represent the capacity of an unreliable machine, in wireless
communication, in which a Markov chain is used to depict randomly
time varying signals or channels. In financial engineering, a
geometric Brownian motion model for a stock is frequently used.
The traditional setup can be described by a linear stochastic
differential equation, where both the appreciate rate and
volatility are constant. However, it has been recognized that such
a formulation is far from realistic. Very often, there are
additional randomness due to the variation of interest rates and
other random environment factors. For example, the well-known
Markowitz's mean-variance portfolio selection is one of the LQ
control problems. Some recent effort for mean-variance control
problems has been on obtaining optimal portfolio selections when
both the appreciation rate and volatility depend on a Markov
chain. For all of the applications mentioned above,
practical considerations often lead to deal with
unobservable Markov chains.
In many situation, the Markov chain is used to model random environment.
Thus, treat adaptive controls,
stability, and stabilization of such systems will  have significant
impact to many applications.

There have been continued interest in dealing with hybrid systems
under a Markov switching. In \cite{Mariton},  stabilization for
robust controls of jump LQ control problems
was investigated. In \cite{JiC}, both controllability and
stabilizability of jump linear LQ systems were considered.
Stability under random perturbations of Markov chain type can be
traced back to the work \cite{KacK}. This line of work has been
substantially
expanded to diffusion systems in \cite{K,Kushner67}.
Recently, renewed interests have been shown to deal with switching
diffusions; see for example \cite{KZY,Mao99,MYY,ZhuY} among  others.

In the literature, stabilization of
continuous-time, adaptive control systems
with hidden Markov chains were considered in \cite{CainesZ,DufourB}.
In both of these references, averaging criteria were used for the
purpose of stabilization. To be more precise,
adaptive control strategies were developed
in \cite{DufourB} to make both the system and the control
have bounded second moment in the sense
$$\limsup_{t\to\infty} \dE[|X(t)|^2 +|U(t)|^2] <\infty,$$
whereas adaptive controls were obtained in \cite{CainesZ} to have
the second moments of
the averages of both the system and control  bounded
in the sense
$$\limsup_{t\to\infty} {1\over t} \dE \left[\int^t_0 [ |X(s)|^2
+|U(s)|^2]
ds\right] <\infty.$$

In comparisons with the aforementioned references, it is a
worthwhile effort to examine the pathwise stabilization of the
associated  LQ problems under partial observations. First, to be
of any practical use in applications, the system resulting from an
adaptive control law should not allow wild behavior in the sample
paths. Secondly, owing to the use of adaptive control strategies,
known results in stability and stabilization in Markov-modulated
stochastic systems cannot be applied directly. As will be seen in
later section, the feedback adaptive controls render difficulty in
analyzing the underlying systems. Certain functions associated
with the diffusion matrix  in fact grow faster than  normally is
allowed in the standard analysis. When averaged criteria are used,
this kind of difficulty will not show up since by taking
expectation, we can easily average out the Brownian motion term.
However, when pathwise criteria are the used, we can no longer use
the argument based on  using expectations. Thus the
consideration of pathwise stabilization is both practically
necessary and theoretically interesting.

To begin our quest of finding admissible controls that stabilize
the systems almost surely, we answer the question if the
controlled process is regular. By a process being regular we mean
that it does not have finite explosion time with probability one. We
establish regularity under feedback controls under linear growth
conditions for the feedback controls. Then, we develop sufficient
conditions and admissible adaptive controls under which the system
is stabilizable. Moreover, as a by-product, we also establish
positive recurrence of the underlying processes as a corollary of
our stabilization result. For a deterministic system given by a
differential equation, if the solutions are ultimately uniformly
bounded, then it is Lagrange stable. For stochastic systems,
almost sure  boundedness excludes many cases (e.g., any systems
perturbed by a white noise). Thus, in lieu of such a boundedness,
one seeks stability in certain weak sense. So a process
is recurrent if it starts from a point outside
 a compact set, the process will return to the bounded
set with probability one. We say the process
is positive recurrent if the expected
return time is finite.
In fact, positive recurrence
is termed weak stability for diffusion processes
in \cite{Wonham}.
For a practical system, no finite explosion
time is a must. In addition,
starting from a point outside of a bounded set, the
system should be able to return to the set infinitely often with
probability one.
Moreover, the average return time cannot be
infinitely long otherwise the controlled system is useless.
Thus, regularity and recurrence of adaptive
control systems can be viewed
as ``practical'' stability conditions.

The rest of the paper is organised as follows. Section 2 presents
the formulation and preliminaries. Section 3 investigates the
regularity of the underlying process. Our conclusion is that, as
long as the feedback controls have linear growth, the resulting
systems will be regular. Section 4 proceeds with the study of
stabilization.  We conclude the paper with some additional
remarks in Section 5. In order to preserve the flow of presentation,
proofs of a couple of technical results are postponed to two
appendices to facilitate the reading.

\section{Formulation and Preliminary}
\subsection{Problem Setup}
Denote by $(\Omega,\A, \dP)$ a probability space with an
associated nondecreasing family of $\sigma$-algebras $(\F_t)$.
Let $\al(t)$ be a continuous-time Markov chain with a finite state
space $\M=\{1,\ldots,m\}$ and transition rate matrix
$\Pi =(\pi_{ij}) \in \rr^{m\times m}$,
and $W(t)$ be a standard $\rr^n$-valued Brownian motion. In
the above and hereafter, $A'$ denotes the transpose of a matrix $A$,
$|A|= \sqrt{\tr(AA')}$ is the trace
norm of $A$, and $|v| =\sqrt {v'v}$ is the usual Euclidean norm of a vector $v$.

Assume throughout the paper that $W(t)$ and $\al(t)$ are
independent. Let $X(t) \in \rr^{n}$ and $U(t)\in \rr^d$ be the
state and control processes, respectively. For $i \in \M$, $A_i \in
\rr^{n\times n}$ and $B_i \in \rr^{n\times d}$ are matrices
with appropriate dimensions. Our main
interest focuses on the following regime-switching stochastic
system
\beq{sys-o}\barray
\ad d X(t)=[ A_{\al(t)} X(t) +  B_{\al(t)} U(t)] dt + dW(t)
\earray\eeq
with square integrable initial condition $X(0)=x$
As in \cite{CainesZ,DufourB},
denoting the column vector of $\rr^{m}$ of indicator functions by
$$\Phi(t)= ( \ind_{\{\al(t)=1)\}}, \ldots, \ind_{\{\al(t)=m\}})'$$
where $\ind_E$ stands for the usual indicator function of the event $E$,
we may present the dynamics of the Markov chain by
$$
d\Phi(t) = \Pi' \Phi(t) dt + d M(t).
$$
The process $M(t)$ is an $\rr^m$-valued square
integrable martingale with right
continuous
trajectories. The independence
of $\al(t)$ and $W(t)$ implies that of $\Phi(t)$
and $W(t)$. In all the sequel, we also assume that $x$,
$\Phi(t)$, and $W(t)$ are mutually independent.
Consider the quadratic cost criterion
$$
J_T(x,\Phi,U) = \dE_{x,\al}\Bigg[ {1\over 2} \int^T_0[ X'(t) Q_{\al(t)}
X(t) + U'(t) R_{\al(t)} U(t))] dt\Bigg],
$$
where $\dE_{x,\al}$ denotes the expectation with initial conditions
$X(0)=x$, $\al(0)=\al$, and for each
$i\in \M$, $Q_i$ is a
symmetric positive semi-definite
matrix, and $R_i$ is a symmetric positive definite matrix.

One of the main features of the system considered here is that the
Markov chain under consideration is a hidden one. As treated in
\cite{CainesZ,DufourB}, the essence is that we are dealing with a
system (\ref{sys-o}) with unknown mode that switches back and
forth among a finite set at random times. But different from
previous consideration, we wish to establish the regularity of the
process and to find conditions ensuring almost sure stabilization.
The almost sure stabilization poses new challenges and
difficulties since we cannot average out the martingale term by
means of taking expectations. Compared with the aforementioned papers,
different techniques are needed. Here the keystone is to find a suitable
Lyapunov function.

Throughout the paper, the process $X(t)$ is
assumed to be observable, but this is not
the case for the switching process $\al(t)$.
The problem belongs to the category
of
controls with partial observations.
Observing $\al(t)$ through the adaptive control
process with Gaussian white noise brings us to the framework of
the setup of Wonham filtering problems
\cite{Wonham65}.
Denote by $\F^X_t$ the
$\sigma$-algebra generated by $\F^X_t =\sigma\{X(s), s \le t\}$.
For the problem of interest, a control is said to be admissible if
for each $t \ge 0$, $U(t)$ is $\F^X_t$-measurable. We are now
in position to state precisely the problem we wish to study.

\para{Problem statement.}
Under the setup presented so far, we aim to solve the following problem.
 \begin{enumerate}
 \item  We  analyze (\ref{sys-o}) and obtain conditions
 under which the system will
 be regular. Hence, our goal is to propose sufficient conditions ensuring
the process will not have finite explosion time.
We show that, as long as the feedback control $U$ (as a function of
$x$) has linear
growth in $x$, the resulting adaptive control system will be regular.
 \item We design
 admissible adaptive controls
 and provide sufficient conditions that stabilize
the closed-loop system almost surely (a.s.). Loosely, the
sufficient condition ensures that for almost all sample points
$\omega$ (except a null set), the corresponding system will be
stabilizable. The precise definition of almost sure stabilization
will be provided in the next section.
\end{enumerate}

\subsection{Preliminary}
As in \cite{CainesZ,DufourB}, we convert this partially observed
system to a control process with complete observation. It entails to
replace the hidden state $\Phi(t)$ by its estimator, namely the
well-known Wonham filter $\wdh \Phi(t)$.
Using feedback control $U(t)=U(X(t),\wdh \Phi(t))$,
we shall need the following notation
\beq{notation}\barray
\ad \wdh \Phi_i(t)= \dE [\ind_{\{\al(t)=i\}}| \F^X_ {t}],\\
\ad \wdh \Phi(t) =(\wdh \Phi_1(t) ,\ldots, \wdh \Phi_m(t))' \in \rr^{m},\\
\ad C(X(t))= (A_1 X(t)+B_1 U(t),\ldots, A_m X(t)+B_m U(t)) \in \rr^{n\times
m},\\
\ad D(\ph ) = (\diag (\ph) -\ph \ph' ) \ \hbox{ for }\ \ph  \in \rr^{m},\\
\ad \diag(\ph) =\diag (\ph_1, \ldots , \ph_m).
\earray\eeq
Denote also the innovation process by
$$
dV(t)= dX(t)- C(X(t)) \wdh \Phi(t) dt.
$$
Using the above notation, we can rewrite the converted completely observable
system as
\beq{sys}
 d \pmatrix{ X(t) \cr \wdh \Phi(t)}
 = \pmatrix{ C(X(t)) \wdh \Phi(t)\cr
 \Pi' \wdh \Phi(t) } dt + \pmatrix{ I_n \cr  D(\wdh \Phi(t))C(X(t))'}
 d V(t),
\eeq
where
$I_n$ stands for the identity matrix of order $n$.

\begin{rem}\bdd
{\rm  Before proceeding further, we shall make a few remarks.
\begin{itemize}
\item
The form $C(X(t))$ indicates the $X(t)$-dependence.
When the feedback control $U(t)$ is of linear form,
it depends on $X(t)$ linearly. This point will be used in what follows.
\item
The equivalent and completely observable system can be viewed as a
controlled diffusion, in which the usual diffusion term is replaced by
$$\pmatrix{ I_n \cr D(\wdh \Phi(t))C(X(t))'}$$
and the driven Brownian motion is given by $V(t)$.
\item
When linear feedback control is used, both the drift and diffusion
grow at most linearly, which is a useful observation.
\item
Since $\wdh \Phi(t)$ is the probability conditioned on the observation,
for each $t\ge 0$ and each $i \in \M$, $\wdh \Phi_i (t) \ge 0$ with
$$ \sum^m_{i=1}\wdh\Phi_i(t)=1.$$
\end{itemize}
}\end{rem}

Denote the joint vector
by $Y(t)=(X(t),\wdh \Phi(t))' \in  \rr^{n+m}$.
In what follows, we often consider
$|Y(t)| \ge r$ for some $r >0$, where $|Y|$ is the
usual Euclidean norm. Denote by $N(0;r) \in \rr^{n+m}$ the
neighborhood centered at 0 with radius $r$. Using the notation
defined in (\ref{notation}) associated with the stochastic
differential equation
(\ref{sys}), we define the following operator. For each
sufficiently smooth real-valued function
$ h: \rr^{n+m}\slash N(0;r) \mapsto \rr$, define
\beq{op-def}\barray \op h(y)\ad=  \op h(x,\ph)\\
\ad =\Bigl(\nabla h(x,\ph)\Bigr)'
\pmatrix{ C(x) \ph \cr
 \Pi' \ph }
  \\  \aad \
 + {1\over 2}\tr\bigg( \pmatrix{ I_n & C(x)D(\ph)'}  \nabla^2 h(x,\ph)
\pmatrix{ I_n & C(x)D(\ph)'}' \bigg)\earray\eeq
where $\nabla h$ and $\nabla^2 h$ are the gradient and Hessian of $h$,
respectively.

\section{Regularity}
First, let us recall the definition of regularity.
According to \cite{K}, the  Markov process
$$Y(t)=\pmatrix{X(t) \cr\wdh \Phi(t)}$$
is regular, if for any $0<T<\infty$,
$$
\dP\left(\sup_{0\le t\le T}|Y(t)|=\infty \right)=0.
$$
Roughly speaking, regularity ensures the process
under consideration will not have finite explosion time.
For our adaptive control systems, we proceed to show that under linear
feedback control, the systems is regular.

\begin{thm}\bdd\label{regularity}
Assume that the feedback control $U(t)=U(X(t),\wdh \Phi(t))$ is admissible and
that
it grows at most linearly in $X(t)$. Then, the feedback
control system {\rm(\ref{sys})} is regular.
\end{thm}

\begin{rem}\bdd\label{linear-g}
{\rm In fact, for our problem, we are mainly interested in linear
(in $x$ variable)
feedback controls. In this case, the linear growth condition
is clearly satisfied.
}\end{rem}

\para{Proof.}
Let $\Theta$ be an open set in $\rr^{n+m}$ and denote
$$O= \Bigl\{ y =(x,\ph)' \in \Theta,  \ph=(\ph_1,\ldots,\ph_m) \hbox{
satisfying } \ph_i \ge 0 \hbox{ for } i \in \M, \hbox{ and }
\sum^m_{i=1} \ph_i =1\Bigr\}.$$
We first observe that both the drift and the diffusion coefficient  given in
(\ref{sys}) satisfy the linear growth and Lipschitz condition
in every open set in $ O \subset \rr^{n+m}$. Thus, to prove the
regularity, using the result in \cite{K}, we only need to show that
there is a nonnegative function $\mathcal{U}$ which is twice
continuously differentiable in $O_{r}= \{y \in O, |y| >r\}$
for some $r>0$ with $y=(x,\ph)'$ such that
\beq{r-cond-1}
\inf_{|y|>R} \mathcal{U}(y) \to \infty \ \hbox{ as } \ R \to \infty,
\eeq
and that there is an $\ga>0$ satisfying
\beq{r-cond-2}  \op \mathcal{U}(y)\le
\ga \mathcal{U}(y).
\eeq
Thus, to verify the regularity of  the process $Y(t)$,
all needed is to construct an appropriate Lyapunov function $\mathcal{U}$.
Note that we only need a Lyapunov function that is smooth and defined in the
complement of a sphere.
Equivalently, we only need the smoothness of the Lyapunov function
to be in a deleted neighborhood of the origin. To this end, take
$r=1$ and denote by $O_1$ the set
\beq{o1-def} \barray O_1\ad =
\Bigl\{ y =(x,\ph)' \in \rr^{n+m},  \ |y| >1  \hbox{ and }
\ph=(\ph_1,\ldots,\ph_m) \hbox{
satisfying } \\
\aad \qquad \qquad  \ph_i \ge 0 \hbox{ for } i \in \M, \hbox{ and }
\sum^m_{i=1} \ph_i =1\Bigr\}.\earray\eeq
Define $\mathcal{U}: O_1\mapsto \rr$ as
$\mathcal{U}(y)=|y|$.  It is easily checked that
condition (\ref{r-cond-1}) holds. Moreover, we have
$$\nabla  \Big| \pmatrix{ x\cr\ph }\Big| = { \pmatrix{ x\cr\ph } \over  \Big|
\pmatrix{ x\cr\ph }\Big|},$$
and
$$\nabla^2  \Big| \pmatrix{ x\cr\ph }\Big|={I_{n+m} \over  \Big|
\pmatrix{ x\cr\ph }\Big|} -  {\pmatrix{ xx'& x \ph'\cr  \ph x' & \ph \ph' }\over
\Big| \pmatrix{ x\cr\ph } \Big|^3}.$$
Consequently, it follows from (\ref{op-def}) that
\beq{l-U}
\barray
\op \mathcal{U}(x,\ph) \ad = {1\over \Big| \pmatrix{ x\cr\ph }
\Big|} (x' C(x) \ph + \ph'
 \Pi' \ph ) \\
 \aad \ +  {1\over 2 \Big| \pmatrix{ x\cr\ph }
\Big|}(n+ \tr (D(\ph)  C'(x) C(x) D'(\ph) ))   \\
\aad \ - {1\over 2\Big| \pmatrix{ x\cr\ph } \Big|^3}\tr\bigg( \pmatrix{ I_n &
C(x)D(\ph)'}  \pmatrix{ xx'& x \ph'\cr  \ph x' & \ph \ph' } \pmatrix{ I_n &
C(x)D(\ph)'}' \bigg).\earray\eeq
Note that the set that we are
working with is $O_1$ defined in (\ref{o1-def}).
In particular,  the use of $O_1$ yields that
for any $y \in O_1$, $|\ph|$ is always bounded.
We also note that owing to the definition of $C(x)$ and the linear
growth feedback controls used,
$C(x)$ is a function grows at mostly  linearly in $x$.
To proceed, henceforth,
use $\gamma$ as a generic positive constant with the
convention that $\gamma+\gamma=\gamma$ and $\gamma\gamma=\gamma$ in an
appropriate sense. It
follows that for the terms on the third line from bottom of
(\ref{l-U}),  for $|(x,\ph)'|$ large enough,
\bea\ad  {1\over \Big| \pmatrix{ x\cr\ph }\Big|} \Bigl|
x' C(x)\ph +
\ph' \Pi' \ph \Bigr|
\le {\gamma \over \Big| \pmatrix{ x\cr\ph }\Big|}
\Big| \pmatrix{ x\cr\ph }\Big|^2
\le \gamma \Big| \pmatrix{ x\cr\ph }\Big|.
\eea
Likewise, for the next two term, we have
\bea \ad {1\over 2 \Big| \pmatrix{ x\cr\ph }
\Big|}\Bigl| n+ \tr (D(\ph)  C'(x) C(x) D'(\ph) ) \Bigr|
\le \gamma \Big| \pmatrix{ x\cr\ph }\Big|.
\eea
Combining the above estimates, we can deduce that
$$\op \mathcal{U}(x,\ph) \le \gamma \Big| \pmatrix{ x \cr \ph}\Big| = \gamma
\mathcal{U}(x,\ph)$$
for some $\gamma>0$. Consequently, the second condition (\ref{r-cond-2}) is
satisfied. Thus
the regularity of the feedback control is obtained.
\qed

\section{Stabilization}
In this section, we establish conditions under which the system of
interest is stabilizable in the almost sure sense. We first
present the definition and then proceed to find
sufficient conditions for stabilization.

\begin{defn}\bdd\label{stab-l}
{\rm System (\ref{sys-o}) or equivalently
(\ref{sys})  is said to be almost surely stabilizable if there is
a feedback control law $U(t)$ such that the resulting trajectories satisfy
\beq{as-s} \limsup_{t\to \infty} {1\over t} \log |X(t)|
\le 0 \ \hbox{ almost
surely.}\eeq
} \end{defn}

Note that the definition given in (\ref{as-s}) is natural. When
studying stability of stochastic differential equations, especially
for pathwise stability, one uses the so-called
$q$th-moment Lyapunov exponent
$$ \limsup_{t\to \infty} {1\over t} \log |X(t)|^q$$
for some $q>0$. Here, roughly, we require that under the control law,
the first-moment Lyapunov exponent is non-positive.

\subsection{Auxiliary Results}
Before proceeding further, let us first recall a lemma, which is
concerned with the existence of the associated system of Riccati
equations when quadratic cost criteria are used. The proof of the
lemma is given in \cite{KFJ}.

\begin{lem}\bdd\label{Riccati}
Consider the system of Riccati equations
\beq{ric} A'_i P_i +P_i A_i -P_i B_i R^{-1} B'_i P_i +
\sum^m_{j=1} \pi_{ij} P_j + Q =0, \ i\in\M,\eeq
  where $Q  \in \rr^{n\times n}$ is symmetric and  positive semi-definite,  and
  $R \in \rr^{m\times m}$ is
symmetric and  positive definite.
The system {\rm (\ref{ric})} has a solution if and only if
for each $i\in \M$,
 there is a matrix $\lbar P_i$ satisfying
\beq{nec-suf} A'_i \lbar P_i + \lbar P_i A_i - \lbar P_i A_i - \lbar P_i
 B_i R^{-1} B'_i \lbar P_i + \sum^m_{j=1} \pi_{ij} \lbar P_j + Q
 \le 0.\eeq
 Furthermore, if $Q$ is positive definite, so are $P_i$ for $i\in
 \M$.
\end{lem}

To carry out the analysis, we need  some auxiliary results on the
bounds of the quadratic variation process.
Before getting the almost sure bounds, we examine
the moment bounds for certain  related martingales, which turn out to be
interesting in their own right.
 The main ingredient is the use of
properties of the associated Markov chain.

\subsubsection*{Moment Bounds}

\begin{prop}\bdd\label{2nd-m}
Consider the stochastic differential equation
\beq{sde-2} d \wdh \Phi(t) = \Pi' \wdh \Phi(t) dt +
D(\wdh \Phi(t)) C(X(t))' dV(t)
\eeq
and define the associate martingale
\beq{mg-2} N(t) = \int^t_0 D(\wdh \Phi(s))
C(X(s))' dV(s).
\eeq
Suppose that the Markov chain $\al(t)$  is irreducible.
Then, for some positive constant $K$ independent
of $t$,
\beq{moment-bd}
\dE \left[{1 \over t}|N(t)|^2\right]\le K.
\eeq
\end{prop}

\para{Proof.} The proof is given in Appendix A. \qed

\begin{rem}\bdd\label{rem-clt}
{\rm It follows from the proof of Proposition~\ref{2nd-m} that the limit of the
matrix
$$S=\lim_{t\to \infty}
{1\over t} \int^t_0 \int^t_0 \Pi' [\wdh \Phi(u)-\nu]
[\wdh \Phi(s)-\nu]' \Pi  du ds $$
is finite. Clearly, this matrix is symmetric and positive semi definite.
A moment of reflect reveals that we can further prove the asymptotic
normality. That is
$${1\over \sqrt t} \int^t_0
\Pi' [\wdh \Phi(s)-\nu] ds \
\hbox{converges in distribution to } \ \mathcal{N}(0,S) \ \hbox{ as }
t\to \infty.
$$
That is, a normalized sequence defined on the left-hand side above
converges in distribution to a normal random vector with mean 0
and covariance $S$.

Another ramification is that in lieu of considering the second-moment
bounds, we can deal with $q$th-moment bounds.
In fact, using the same techniques, we can show that
for any integer $p>0$,
$$ \dE\left[ \l {1\over \sqrt t}\int^t_0 \Pi' [\wdh \Phi (s) -\nu] ds
\r^{2p}\right] <
\infty.$$
Hence, as the solution (\ref{sde-2}) is given by
$$
\wdh \Phi(t) = \wdh \Phi(0) + \int^t_0 \Pi' \wdh \Phi (s) ds + N(t)
$$
which means that
$$
N(t)=\wdh \Phi(t)- \wdh \Phi(0)-\int^t_0 \Pi' \wdh \Phi (s) ds,
$$
we obtain that
$$
\dE\left[ \l{1\over \sqrt t} N(t) \r^{2p} \right]< \infty.$$
Next, for odd exponents and for any integer $p\geq 1$, it follows from
H\"older's inequality that
$$ \left(\dE \left[\l{1\over\sqrt t} N(t) \r^{2p-1}\right]\right)^{2p}\le
\left( \dE \left[ \l{1\over\sqrt t} N(t) \r^{2p}\right]\right)^{2p -1}
 < \infty.$$
Finally, we conclude that for any positive integer $q$,
$$ \dE\left[ \l{1\over\sqrt t} N(t) \r^q\right]<\infty.$$
}
\end{rem}

\subsubsection*{Almost Sure Bounds}

For the almost sure stabilization, we need to show that
$${1\over t}  |N(t)|^2   \le K \hspace{1cm} \hbox{ a.s.}$$
for some $K>0$ independent of $t$.

\begin{prop}\bdd\label{mg-2-c}
Consider {\rm (\ref{sde-2})} and
suppose that the Markov chain $\al(t)$  is irreducible.
Then, the quadratic variation of the process
$N(t)$ satisfies $\lan N, N\ran_{t} \le K t$ where $K$ is some positive constant
independent of $t$.
Therefore,
\beq{mg-0} \lim_{t \rightarrow \infty}{1\over t} N(t)=0 \hspace{1cm} \hbox{ a.s.}
\eeq
\end{prop}

\para{Proof.} The proof is given in Appendix B. \qed

\subsection{Stabilization}

\begin{lem}\bdd\label{lv}
Consider the set $\Delta$ defined by
$$
\barray
\Delta \ad = \Bigl\{ (x,\ph) \in \rr^{n}\times \rr^{m},
\ph=(\ph_1,\ldots,\ph_m) \hbox{ satisfying } \ph_i \ge 0 \hbox{ and }
\sum^m_{i=1} \ph_i =1\Bigr\}.
\earray
$$
Denote
\beq{ph-def-i} P(\ph)= \sum^m_{i=1} P_i \ph_i.\eeq
For some $\theta >0$, let $V_{\theta}(x,\ph): \Delta \mapsto \rr$ with
$V_{\theta}(x,\ph)= \log  (\theta + x' P( \ph) x) $.
Then, we have
\beq{d-log}
\barray
\op V_{\theta}(x,\ph) \ad =
{1\over \theta+ x'P(\ph) x } \bigg(2x'P(\ph)C(x)\ph +  (x' \wdt P x)'  \Pi'
\ph\bigg)\\
\aad \ + {1\over \theta+ x'P(\ph) x }
\hbox{\rm tr} \bigg(P(\ph)+2C(x)D(\ph)'x'\wdt P\bigg)    \\
\aad \ - {1\over 2 (\theta+x'P(\ph) x )^2}
\hbox{\rm tr} \bigg( \bigl( I_{n} \: , \: C(x) D(\ph)' \bigr) \Lambda(x,\ph)
\bigl( I_{n}
\: , \: C(x) D(\ph)' \bigr)' \Bigg),
\earray\eeq
where\\
$$\Lambda(x,\ph)=\pmatrix{ 2 P(\ph) x \cr x' \wdt P x } \pmatrix{ 2 P(\ph) x \cr
x' \wdt P x }',$$
$$P= (P_1,\ldots,P_m)', \quad x'\wdt P x = (x'P_1 x,\ldots, x' P_m x)' \in
\rr^{m},$$
$$Px = (P_1 x,\ldots,P_m x) \in \rr^{n\times m}, \quad x' \wdt P= (\wdt P x)'
\in \rr^{m\times n}.$$
\end{lem}
\para{Proof.}
We have
$$ \nabla \log  (\theta + x' P( \ph) x) =
{\pmatrix{ 2P(\ph) x \cr x' \wdt P x }\over  \theta+x'P(\ph)x  }
$$
and
$$\nabla^2 \log  (\theta + x' P( \ph) x)
= - {
 \pmatrix{ 2 P(\ph) x \cr x' \wdt P x }
  \pmatrix{ 2 P(\ph) x \cr x' \wdt P x }' \over (\theta + x' P(\ph) x)^2}
+ { 2  \pmatrix{ P(\ph)  &  \wdt P x  \cr x' \wdt P
& 0_{m} }   \over \theta + x' P(\ph) x}
\vspace{2ex}
$$
where $0_{m}$ stands for a square matrix of order $m$ with all entries equal to
zero. Consequently, it follows from (\ref{op-def}) that
\beq{}\barray
\op V_{\theta}(x,\ph) \ad = {1\over    \theta + x'P(\ph) x } \Bigl( 2x' P(\ph)
C(x) \ph  +  (x' \wdt P x)'  \Pi' \ph \Bigr)  \\
\aad \quad + {1\over 2}
\hbox{\rm tr} \Bigl( \bigl( I_{n} \: , \: C(x) D(\ph)' \bigr) \ \nabla^2
V_{\theta}(x,\ph)
\ \bigl( I_{n} \: , \: C(x) D(\ph)' \bigr)' \Bigr),
\earray\eeq
which immediately implies (\ref{d-log}). \qed

\noindent
For the purpose of stabilization, we also need an estimate on $\op
V_{\theta}(X(t), \wdh \Phi(t))$.

\begin{lem}\bdd\label{lv<0}
Assume that equation {\rm(\ref{nec-suf})} is satisfied and that
$$Q - {1\over 2} \Bigl[P_{i}B_{i}-P_{j}B_{j}\Bigr] R^{-1}\Bigl[P_{i}B_{i}-
P_{j}B_{j}\Bigr]'$$
are positive definite matrices for all $(i,j)\in\M^{2}$ where $P_i$ for $i\in\M$
are the
solutions of the algebraic Riccati equations given by {\rm(\ref{Riccati})}.
Then, the infinitesimal generator of the process $(X(t),\wdh \Phi(t))$
associated with the feedback control law
\beq{u-cont}
U(t) = - R^{-1} \sum^m_{i=1} \wdh \Phi_i(t) B'_i P_i X(t),
\eeq
satisfies for some constant $\gamma>0$
\beq{majlv}
\op V_{\theta}(X(t), \wdh \Phi(t)) \le {\gamma \over \theta}.
\eeq
\end{lem}

\para{Proof.} We can deduce from Lemma \ref{lv} that
\begin{eqnarray*}
\barray
\op V_{\theta}(X(t), \wdh \Phi(t)) \ad
\leq
{1\over \theta+ X(t)'P(\wdh \Phi(t)) X(t) } (2X(t)'P(\wdh \Phi(t))
C(X(t)) \wdh \Phi(t) )\\
\aad\  + {1\over \theta+ X(t)'P(\wdh \Phi(t)) X(t) }
((X(t)' \wdt P X(t))'  \Pi' \wdh \Phi(t))\\
\aad\ +{1\over \theta+X(t)'P(\wdh \Phi(t)) X(t)}
\tr \bigg(P(\wdh \Phi(t))+2C(X(t))D(\wdh \Phi(t))'X(t)'\wdt P\bigg).
\earray\end{eqnarray*}
Therefore, following exactly the same lines as in \cite{DufourB},
we obtain that
\begin{eqnarray*}
\barray
\op V_{\theta}(X(t), \wdh \Phi(t)) \ad
\leq
-{ 1 \over \theta + X(t)'P(\wdh \Phi(t)) X(t)}
\Bigl( X(t)'
\Bigl[ Q \\
\aad\ - \sum^m_{i=1}\sum^m_{j=1} { \wdh \Phi_i(t) \wdh \Phi_j(t)
\over 2}
\bigg[P_iB_i -P_jB_j\bigg]
 R^{-1}\bigg[P_iB_i-P_jB_j\bigg]'\\
\aad\ + \Bigl( \sum^m_{j=1} \wdh \Phi_j(t) B'_j P_j\Bigr)' R^{-
1}\Bigl(\sum^m_{i=1}
\wdh \Phi_i(t)  B'_i
P_i\Bigr) \Bigr] X(t) - \tr (P(\wdh \Phi(t))) \Bigr).
\earray\end{eqnarray*}
Finally,
$$\op V_{\theta}(X(t), \wdh \Phi(t)) \le {1\over \theta} \sum_{i=1}^{m} \tr
(P_i)$$
which completes the proof of Lemma~\ref{lv<0}.\qed

\begin{thm}\bdd\label{stab}
Assume that the conditions of
\lemref{lv<0}
are satisfied. Then, the feedback
control law defined in equation {\rm(\ref{u-cont})} stabilizes the system {\rm
(\ref{sys})} almost surely.
\end{thm}

\para{Proof.}
It follows from Ito's rule that
\beq{ito} \barray
V_{\theta}(X(t),\wdh \Phi(t))\ad= V_{\theta}(x, \ph)+ \int^t_0 \op
V_{\theta}(X(s),\wdh \Phi(s)) ds + M(t)
\earray\eeq
with the initial condition $X(0)=x$ and $\wdh \Phi(0)=\ph$ and the martingale
term
$$M(t)=\int^t_0 \Sigma(s) d V(s)$$
where
 \bea \Sigma(s)\ad = {1\over   \Theta + X(s)'P(\wdh \Phi(s)) X(s) }
   \pmatrix{ 2 X(s)' P(\wdh \Phi(s))  & (X(s)' \wdt P X(s))' }
\pmatrix{ I_n \cr D(\wdh \Phi(s))C(X(s))'},\\
\ad ={1\over    \theta + X(s)'P(\wdh \Phi(s)) X(s) }
\bigl(2X(s)'P(\wdh \Phi(s)) +(X(s)' \wdt P X(s))' D(\wdh \Phi(s))C(X(s))' \bigr).
\eea
We can split the martingale $M(t)$ into two terms,
$M(t)=N_1(t)+N_2(t)$ with
\begin{eqnarray*}
N_1(t)&=&\int^t_0
{2X(s)'P(\wdh \Phi(s))  \over
\theta + X(s)'P(\wdh \Phi(s)) X(s) }dV(s),\\
N_2(t)&=& \int^t_0
{(X(s)' \wdt P X(s))'   \over
\theta + X(s)'P(\wdh \Phi(s)) X(s) }D(\wdh \Phi(s))C(X(s))' dV(s).
\end{eqnarray*}
It is easy to see that
$$ {4X(t)'P(\wdh \Phi(t))P(\wdh \Phi(t)) X(t)  \over (\theta + X(t)'P(\wdh
\Phi(t)) X(t))^{2} }  \le
K_1
\hspace{0.5cm}\textrm{where}\hspace{0.5cm}
K_1={m \over \theta}\max_{i \in \M} (\lambda_{max}(P_i)).
$$
Then, the quadratic variation of $N_1(t)$ satisfies
$\lan N_{1},N_{1} \ran_{t} \le K_1t$ a.s.
Consequently, we deduce from the strong law of large numbers for local
martingales
\cite{Lip} that
\beq{lim-1}
 \lim_{t \rightarrow \infty}{1\over t} N_1(t)=0 \hspace{1cm} \hbox{ a.s.}
\eeq
In view of \propref{mg-2-c}, one can also find a positive constant $K_2$,
independent of $t$, such that
\beq{est-1}\barray
\lan N_2, N_2\ran_{t} \ad = \int^t_0 {|X(s)' \wdt P X(s)|^2  \over
(\theta + X(s)'P(\wdh \Phi(s)) X(s))^2 }|D(\wdh \Phi(s))C(X(s))'|^2 ds\\
\ad \le K_2 t \hspace{1cm} \hbox{ a.s.} \earray\eeq
It also ensures that
\beq{lim-2}
 \lim_{t \rightarrow \infty}{1\over t} N_2(t)=0 \hspace{1cm} \hbox{ a.s.}
\eeq
Therefore, (\ref{lim-1}) and (\ref{lim-2}) imply that
\beq{lim-3}
 \lim_{t \rightarrow \infty}{1\over t} M(t)=0 \hspace{1cm} \hbox{ a.s.}
\eeq
Thus, we find from (\ref{ito}) that
$${1\over t} V_{\theta}(X(t),\wdh \Phi(t)) = {1\over t} V_{\theta}(x,\ph)
+ {1\over t} \int^t_0 \op V_{\theta}(X(s),\wdh \Phi(s)) ds + o(1) \hspace{1cm}
\hbox{ a.s.}$$
Moreover,
$V_{\theta}(x,\ph)/t=o(1)$ as $t\to \infty$
a.s.
By virtue of \lemref{lv<0}, it follows that for all $\theta >0$
\beq{lim-sup1}
 \limsup_{t\to\infty}
{1\over t} V_{\theta}(X(t),\wdh \Phi(t))  = \limsup_{t\to\infty}
{1\over t} \int^t_0 \op V_{\theta}(X(s),\wdh \Phi(s)) ds
\le {\gamma \over \theta} \hspace{1cm} \hbox{ a.s.}
\eeq
Furthermore, one can observe that
$x'P(\ph) x \ge \lambda_{\min}(P(\ph)) |x|^2$ and since $P(\ph)$ is positive
definite, the
minimal eigenvalue of $P(\ph)$ is positive.
Consequently,
$$\log(\lambda_{\min}(P(\ph)))+ 2\log(|x|) \le
\log(\theta+ \lambda_{\min}(P(\ph)) |x|^2)\le \log(\theta +x'P(\ph) x )$$
which leads to
\beq{lim-sup2}
{1\over t}\Bigl( \log (\lambda_{\min}(P(\wdh \Phi(t)))) +2 \log( |X(t)|) \Bigr)
\le {1\over t} V_{\theta}(X(t),\wdh \Phi(t)).\eeq
Finally, we conclude from (\ref{lim-sup1}) and (\ref{lim-sup2})
that for all $\theta >0 $,
$$\limsup  {1\over t} \log |X(t)| \le {\gamma \over 2 \theta} \hspace{1cm}
\hbox{ a.s.}$$
We complete the proof of Theorem~\ref{stab} by taking the limit as $\theta$
tends to infinity. \qed

\begin{rem}\bdd\label{about-a-est}
{\rm Normally, dealing with stochastic differential equations, to
obtain the almost sure bounds of the solutions, one often relies
on the use of appropriate Lyapunov functions to have the diffusion
term of the process be bounded after a transformation. Here, we
are dealing with a martingale term with some what faster rate of growth in
$x$. Nevertheless, thanks to the second component of the diffusion
(\ref{sde-2}), the probabilistic meaning of $\wdh \Phi(t)$ enables
us to work around the obstacle. To obtain the desired bounds, an
alternative is to obtain an almost sure central limit theorem.
Here however, we take a different approach. The main point is the
use of \propref{mg-2}.
}\end{rem}

Recall the notion of recurrence for the diffusion
process $(X(t),\wdh \Phi(t))$ starting at $X(0)=x$ and $\wdh
\Phi(0)=\ph$. Consider an open set $O$ with compact closure, and let
$$\sigma^{x,\ph}_O=\inf\Bigl\{t>0,(X(t),\wdh \Phi(t))\in O\Bigr\}$$
be the first entrance time of the
diffusion to the set $O$. If $(X(t),\wdh \Phi(t))$ is regular, it
is recurrent with respect to $O$ if $\dP\{\sigma^{x,\ph}_O<\infty\}=1$
for any $ (x,\ph) \in O^c $,  where $O^c$ is the
complement of $O$.
A recurrent process with finite mean recurrence time for some set $O$,
is said to be positive recurrent with respect to $O$, otherwise,
the process is  null recurrent with respect to $O$.
It has been proven in \cite{K} that recurrence and
positive recurrence are independent of the set $O$ chosen. Thus, if
it is recurrent (resp. positive recurrent) with respect to $D$, then it is
recurrent (resp. positive recurrent) with respect to any other open set
$\Theta$ in the domain of interest.
Looking over the proof of the stabilization presented, we could show that
for the Lyapunov function
$$V_{0}(x,\ph)= \log ( x'P(\ph) x),$$
one can find $\gamma>0$ such that for all $ (x,\ph) \in O^c $,
\beq{lp-v<0}\op
V_{0}(x,\ph)
\le -\gamma.\eeq
In view of the known result of positive
recurrence of diffusion processes
\cite{K},  (\ref{lp-v<0}) is precisely a
necessary and sufficient condition
for positive recurrence.
Thus, we obtain the following result as
a by-product.

\begin{cor}\label{p-recu}
Under the conditions of \thmref{stab},
with the control law {\rm (\ref{u-cont})} used, the diffusion systems
{\rm (\ref{sys})} is positive recurrent.
\end{cor}

We would like to add that the positive recurrence of the process
is an important property. It has engineering implication for
various applications. Essentially, it ensures that starting from a
point outside of a bounded set, the control laws enables the
system to return to a compact set almost surely. This may be
viewed as a practical stability condition. In fact, Wonham used
the term weak stability for such a property in his paper
\cite{Wonham}.

\section{Further Remarks}
This paper has been concerned with stabilization in the almost sure
sense of an adaptive control system with linear dynamics modulated by an
unknown Markov chain. Under the framework of Wonham filtering, the
underlying system is converted to a fully observable system. Using
feedback control that is linear in the continuous state variable,
we establish pathwise stabilization of the process.
Along the way of our study, we have also obtained regularity of the
underlying process. In addition,
as a corollary, we have shown that under the
stabilizing control law, the resulting system is positive recurrent.
These results pave a way for practical consideration of stabilization
of adaptive controls of LQ systems with a hidden Markov chain.
Several directions may be worthwhile for further study and
investigation.

\begin{itemize}
\item In our study, irreducibility of the Markov chain is used. We
note that the irreducibility  ensures the spectrum gap
condition or exponential decay in (\ref{mix-1}) and (\ref{mix-2})
of \propref{2nd-m} to hold. It will be interesting to see if
it is possible to remove
 this condition. Our initial thoughts are:
 Under certain conditions, this might be possible.
  For example, if the Markov chain has several irreducible classes
  such that the states in each class vary rapidly, and among
  different classes, they change slowly.
  One may be able to use the
  different time scales to overcome the difficulty under the
  framework of time-scale separation using a
  singular perturbation approach. However, the details on this
  need to be thoroughly worked out; they are in fact
  out of the scope of the current paper.

\item
It will be interesting to design admissible controls
and find sufficient conditions for stabilization of  LQ systems with a
hidden Markov in discrete-time.

\item
In our setup,
the process $X(t)$ represents the noisy observation--hidden Markov
chain observed in white noise. A class of  controlled
regime-switching diffusion systems provides a somewhat more
complex setup. In such a system, the dynamics are represented by
switching diffusions with a hidden Markov chain. The Markov chain
is not observable but can only be observed in another Gaussian
white noise. That is, let us consider the controlled system
\beq{n-sys}\barray
\ad dY(t)= [A_{\al(t)} Y(t) +B_{\al(t)} U(t)]dt +
\sigma_{\al(t)} d V(t)\\
\ad dX(t) = g_{\al(t)} dt + \rho(t) dW(t),
\earray\eeq
where $Y(t)$ and $X(t)$ are vector-valued processes with compatible
dimensions representing the state and observations,
respectively, $V(t)$ and $W(t)$ are independent multi-dimensional
Brownian motions, and $\al(t)$ is the hidden Markov chain with a
finite state space.
As was alluded to in the introduction, one of the
motivations is Markowitz's mean-variance portfolio selections
\cite{ZhouY}.
One may then pose similar stabilization
problems.

\item
Recently, using regime-switching jump diffusions,
which are switching diffusions
with additional external jumps  of
a compound Poisson process, for modeling surplus in insurance risk
has drawn much attention. A related problem in the adaptive setup
is a regime-switching jump diffusion system in which the hidden Markov
chain is observed similar to
the observation in (\ref{n-sys}).
One may then proceed with the study of stabilization problems.

\item
In the study of
stabilization, positive definiteness of certain matrices is
used (see \lemref{Riccati}). A challenging problem is to
investigate the stabilization problem with the positive
definiteness  removed for the system given by (\ref{n-sys}).
Here, the crucial point seems to rely on recent
developments in LQ problems
 with indefinite control weights \cite{CLZ}. One needs to use the
 backward stochastic differential equations from the toolbox of
 stochastic analysis.
\end{itemize}

All of these problems deserve further study and
investigation.


\section*{Appendix A.}

\renewcommand{\thesection}{\Alph{section}}
\renewcommand{\theequation}
{\thesection.\arabic{equation}} \setcounter{section}{1}
\setcounter{section}{1}
\setcounter{equation}{0}

This appendix is devoted to the proof of Proposition~\ref{2nd-m}.
It is divided into several steps.

\noindent
{\bf Step 1.} We already saw that the solution (\ref{sde-2}) is given by
$$
\wdh \Phi(t) = \wdh \Phi(0) + \int^t_0 \Pi' \wdh \Phi (s) ds + N(t).
$$
Consequently
\beq{solu-sde2a}
N(t)=\wdh \Phi(t)- \wdh \Phi(0)-\int^t_0 \Pi' \wdh \Phi (s) ds.
\eeq
In view of (\ref{solu-sde2a}), the probabilistic interpretation of $\wdh
\Phi(t)$ implies that $N(t)$ is a martingale bounded almost surely
for each $t>0$. We proceed to obtain the
moment bounds of $N(t)$.

\noindent
{\bf Step 2.}
As $\Pi$ is the generator of the irreducible Markov chain $\al(t)$,
its unique stationary distribution $\nu$ satisfies
$ \Pi' \nu=0$. Hence, it follows that
$$\int^t_0 \Pi' \wdh \Phi(s)
ds =  \int^t_0 \Pi'  (\wdh \Phi(s) -\nu)
ds.$$
On the one hand, we clearly have from
(\ref{mg-2})
\beq{Esq}
\dE[|N(t)|^2] =\dE\left[ \int^t_0 |D(\wdh \Phi(s)) C(X(s))'|^2 ds \right].\eeq
On the other hand, we deduce from (\ref{solu-sde2a}) that
\beq{mt-2}
\barray\disp {1\over t}\dE[ |N(t)|^2]\ad={1\over t}
\dE\left[ \Bigl| \wdh \Phi(t)- \wdh \Phi(0)- \int^t_0 \Pi' (\wdh \Phi(s)-\nu)
ds \Bigr|^2\right],\\
\ad \le {2\over t} \dE\Bigl[| \wdh \Phi(t)-\wdh \Phi(0)|^2\Bigr] +{2\over t}
\dE \left[ \Bigl| \int^t_0 \Pi' (\wdh \Phi(s)-\nu) ds \Bigr|^2\right],\\
\ad \le {2\over t}+{2\over t}E \int^t_0 \int^t_0
\tr \{ \Pi'\Pi (\wdh \Phi (r)-\nu) (\wdh \Phi' (s)-\nu')\} dr ds.
\earray\eeq
Consider the symmetric matrix
$$G(r,s)=(g_{ij}(r,s))= \dE[(\wdh \Phi(r)-\nu)( \wdh \Phi'(s)-\nu')].$$
One can observe that
\beq{lim-1a}\barray
g_{ij}(r,s) \ad= \dE[(\wdh \Phi(r)-\nu)_i (\wdh \Phi' (s)-\nu')_j], \\
\ad  = \dE[ (\dE[\ind_{\{\al(r)=i\}}| \F^X_r]-\nu_i )(\dE
[\ind_{\{\al(s)=j\}}| \F^X_s]-\nu_j)],\\
\ad = \dE[ \dE[\ind_{\{\al(r)=i\}}| \F^X_r] \dE[\ind_{\{\al(s)=j\}}| \F^X_s]] -
 \nu_j\dE[ \dE[\ind_{\{\al(r)=i\}}| \F^X_r]]  \\
 \aad \qquad \ - \ \nu_i \dE[ \dE[\ind_{\{\al(s)=j\}}| \F^X_s]]
+ \nu_i\nu_j,\\
\ad  = \dE[ \dE[\ind_{\{\al(r)=i\}}| \F^X_r] \dE[\ind_{\{\al(s)=j\}}| \F^X_s]] -
\nu_j \dP(\al(r)=i)
\\ \aad \qquad \ - \
\nu_i \dP(\al(s)=j)
+ \nu_i \nu_j. \earray\eeq
Note also by the Fubini Theorem that
\bea \disp {1\over t} \int^t_0\int^t_0 g_{ij}(r,s) dr ds
\ad ={1\over t}\( \int^t_0\int^t_rg_{ij}(r,s) dr ds +
 \int^t_0\int^r_0 g_{ij}(r,s) dr ds\), \\
\ad =
{1\over t} \int^t_0 \( \int_ r^t
g_{ij}(r,s) ds\) dr
+
{1\over t} \int^t_0\( \int_ s^t g_{ij}(r,s) dr \) ds,\\
\ad = g_1(t)+ g_2(t)=2 g_1(t).
\eea
We have the decomposition
$$ g_1(t)=h_1(t)+\ell_1(t) $$ where
\bea \ad
h_1(t)=
{1\over t} \int^t_0 \( \int^t_r  h(r,s) ds\) dr,
\\
\ad
\ell_1(t)=
{1\over t} \int^t_0 \( \int^t_r  \nu_i (\nu_j - \dP(\al(s)=j))ds\) dr,
\eea
with
\bea
\ad  h(r,s)=
\dP(\al(r)=i)\dP( \al(s)=j | \al(r)=i) - \nu_j P(\al(r)=i).
\eea
Before proceeding further, let us first note the following mixing
properties regarding the Markov chain $\alpha(t)$.
For all $ t\ge 0$ and $s\le t$, denote
\bea \ad p(t)= (\dP(\al(t)=1), \ldots, \dP(\al(t)=m))'\in \rr^m, \\
\ad P(t,s)=( (\dP(\al(t)=j | \al(s)= i),\ i,j\in \M) \in
\rr^{m\times m}, \eea which are the probability vector and
transition matrix of the Markov chain $\alpha(t)$, respectively.
Since $\al(t)$ is irreducible, it is ergodic. Consequently, as $t$
goes to infinity, for the solution of the system \beq{prob-v}
\left\{
\begin{array}{lcl}
{\displaystyle {d p(t) \over dt}}&=& \Pi' p(t) \\
p(0)&=& p_0
\end{array} \right.
\eeq
satisfying
$$p_{0,i} \ge 0 \ \hbox{ and } \sum^m_{i=1}
p_{0,i}=1,$$
one can find two positive constants $\kappa$ and $K$ such that
$p(t) \to \nu$ and
\beq{mix-1}
| p(t)- \nu | \le K \exp(-\kappa t)
\eeq
By virtue of (\ref{mix-1}),
it is easily seen that
\beq{lim-c1}\barray
|\ell_1(t)| \ad =
\l {\nu_i \over t} \int^t_0 \(\int^t_r  (\nu_j - \dP(\al(s)=j)) ds \) dr \r, \\
\ad \le
{\nu_i \over t} \int^t_0 \( \int^t_u  |\nu_j - \dP(\al(s)=j)| ds \) dr, \\
\ad \le {\nu_i K \over t} \int^t_0 \( \int^t_r \exp(-\kappa s) ds \) dr,\\
\ad \le {\nu_i K \over \kappa t} \int^t_0 \exp(-\kappa r) dr,\\
\ad \le {\nu_i K \over \kappa^2 t}.
\earray\eeq
Consequently, $\ell_1(t)$ goes to zero as $t$ tends to infinity.
Next, we shall show that $h_1(t)$ is bounded.
As before, the solution of the system
\beq{prob-m}
\left\{
\begin{array}{lcl}
{\displaystyle {\partial P(t,s) \over \partial t}}&=& \Pi' P(t,s)\\
P(s,s)&=&I_m
\end{array} \right.
\eeq
with $s\le t$, also satisfies for two positive constants
$\lambda$ and $K$, $P(t,s) \to \one\nu'$ and
\beq{mix-2}
| P(t,s)- \one \nu'|
\le K \exp(-\lambda (t-s)).
\eeq
It follows from (\ref{mix-2}) that
\bea
|h_{1}(t)|\ad=
\l {1\over t} \int^t_0 \dP(\al(r)=i) \(\int^t_r
(\dP(\al(s)=j| \al(r)=i) - \nu_j ) ds \)dr \r, \\
\ad
\le
{1\over t} \int^t_0 \(\int^t_r
|\dP(\al(s)=j| \al(r)=i) - \nu_j | ds \)dr, \\
\ad \le {K \over t} \int^t_0 \( \int^t_r \exp(-\lambda(s-r)) ds \) dr,\\
\ad \le {K \over \lambda t} \int^t_0 dr,\\
\ad \le {K \over \lambda}.
\eea
Therefore, $h_1(t)$ as well as $g_1(t)$ are bounded sequences
which ensures that for some positive constant $K$ independent of $t$
\beq{bound-E}
\l \dE \left[
{1\over t} \int^t_0 \int^t_0 \tr \{ \Pi' \Pi (\wdh \Phi(r)-\nu)
(\wdh \Phi'(s) -\nu') \} dr ds  \right]  \r \le K.
\eeq
Finally, (\ref{Esq}) together with (\ref{mt-2}) and (\ref{bound-E}) imply
(\ref{moment-bd})
which completes the proof of Proposition~\ref{2nd-m}.   \qed


\section*{Appendix B.}

\renewcommand{\thesection}{\Alph{section}}
\renewcommand{\theequation}
{\thesection.\arabic{equation}} \setcounter{section}{2}
\setcounter{section}{2}
\setcounter{equation}{0}

We shall now focus on the proof of Proposition~\ref{mg-2-c}. First
of all, we know that
$$\sup_{t\ge 0} | \Pi' \wdh \Phi(t) | \le 1
\hspace{1cm}\hbox{a.s.}$$
In addition, we also have
$|\wdh \Phi(t)|\le 1$ a.s.
Consequently, it follows from
(\ref{solu-sde2a}) that
\beq{m-bd}\barray
| N(t) | \ad \le | \wdh \Phi(t)- \wdh \Phi(0)| + \l\int^t_0 \Pi'
\wdh \Phi(s) ds \r \le 1  + t
\hspace{1cm}\hbox{a.s.}
\earray\eeq
For each $i\in \M$, denote
\bea \ad N_i(t) =  \int^t_0
\sum^n_{j=1} [D(\wdh \Phi(s) C(X(s))']_{ij} d
V_j(s)
\eea
where $[D(\wdh \Phi(s) C(X(s))']_{ij}$ is the $ij$th entry of
the matrix $D(\wdh \Phi(s)) C(X(s))'$ and
$V_j(s)$ stands the $j$th component of $V(s)$.
It follows from the well-known
Doob's martingale inequality given for example in \cite[Theorem
1.7.4, p. 44]{Mao} that for each $i\in \M$ and each positive integer $n$,
\beq{exp-mg}\barray
\ad \dP\( \sup_{0\le t \le n} \l  \int^t_0
N_i(t)-|(D(\wdh \Phi(s)) C(X(s))')_{i,. }|^2 ds   \r \ge
\log n \) \le {1\over n^2},\earray\eeq
where $(D(\wdh \Phi(s) C(X(s))')_{i,. } $ denotes the row vector
in the $i$th row of the matrix $D(\wdh \Phi(s) C(X(s))'$.
Hence, we deduce from the Borel-Cantelli Lemma that for almost all $\omega \in
\Omega$, there is a $K_1=K_1(\omega)>1$ such that for all
$n\ge K_1$ and $t\leq n$
\beq{bd-2a}\barray \disp
\int^t_0
|(D(\wdh \Phi(s) C(X(s))')_{i,.}|^2 ds \ad \le  \log n + N_i(t)
\hspace{1cm}\hbox{a.s.} \\
\ad  \le \log n + 1 + t \hspace{1cm}\hbox{a.s.}
\earray\eeq
The last line above follows from (\ref{m-bd}).
Dividing both sides of (\ref{bd-2a}) by $t$, we obtain that for $n \ge K_2$,
$n-1 \le t\le n$, so
\beq{bd-2b}\barray \disp {1\over t}\int^t_0
|(D(\wdh \Phi(s) C(X(s))')_{i,.}|^2 ds\ad \le {1\over n-1} (\log n + 1 + t)
\hspace{0.5cm}\hbox{a.s.}\\
\ad \le {1\over n-1} (\log n + 1 + n) \hspace{0.5cm}\hbox{a.s.}\\
\ad \le K_3  \hspace{1cm}\hbox{a.s.}\earray\eeq
and the bound $K_{3}$ is independent of $t$.
Consequently, for some positive constant $K$ independent of $t$,
the quadratic variation of the martingale is bounded
by $K t$ almost surely.
That is, (\ref{bd-2b}) implies that
$\lan N , N \ran_{t} \le K t$ a.s.
Finally, we deduce from the strong law of large numbers for local martingales
\cite{Lip} that
$$ \lim_{t \rightarrow \infty}{1\over t} N(t)=0 \hspace{1cm} \hbox{ a.s.}$$
which concludes the proof of Proposition~\ref{mg-2-c}.   \qed



\end{document}